\documentclass[10pt]{amsart}
\usepackage{amsmath,amsthm,amsfonts,amssymb,eucal, hyperref}

\usepackage[T1]{fontenc}
\usepackage[english]{babel}
\usepackage{amsmath,amssymb,amsfonts,amsthm}
\usepackage[pdftex]{graphicx}

\numberwithin{equation}{section}

\newcommand{\de}{\, \mathrm{d}}
\newcommand{\D}{\partial}
\newcommand{\R}{\mathbb R}

\newcommand{\T}{\mathbb T}

\newtheorem{thm}{Theorem}[section]
\newtheorem{lem}[thm]{Lemma}
\newtheorem{prop}[thm]{Proposition}
\newtheorem{cor}[thm]{Corollary}
\newtheorem{rem}[thm]{Remark}
\newtheorem*{open}{Open problem}
\theoremstyle{definition}

\newcommand{\bee}{\begin{equation}}
\newcommand{\ene}{\end{equation}}
\newcommand{\bees}{\begin{equation*}}
\newcommand{\enes}{\end{equation*}}
\newcommand{\bes}{\begin{split}}
\newcommand{\ens}{\end{split}}

\newcommand{\bet}{\begin{thm}}
\newcommand{\ent}{\end{thm}}
\newcommand{\bel}{\begin{lem}}
\newcommand{\enl}{\end{lem}}
\newcommand{\bec}{\begin{cor}}
\newcommand{\enc}{\end{cor}}
\newcommand{\becl}{\begin{cla}}
\newcommand{\encl}{\end{cla}}
\newcommand{\bep}{\begin{proof}}
\newcommand{\enp}{\end{proof}}
\newcommand{\ber}{\begin{rem}}
\newcommand{\enr}{\end{rem}}

\newcommand{\la}{\lambda}

\newcommand{\Z}{\mathbb Z}

\newcommand {\bx}{\mathbf x}

\newcommand {\bk}{\mathbf k}

\newcommand{\bigo}[1]{O\left( #1 \right)}

\begin{document}

\title[Density of states]{A generalised Gauss circle problem and integrated density of states}
\author[Jean Lagac\'e \& Leonid Parnovski]{Jean Lagac\'e \& Leonid Parnovski}
\begin{addresses}
\address{D\'epartement de math\'ematiques et de statistique \\ Universit\'e de Montr\'eal \\ C. P. 6128,
Succ. Centre-ville \\ Montr\'eal, QC\\ H3C 3J7 \\ Canada}
\email{lagacej@dms.umontreal.ca}
\address{Department of Mathematics\\ University College London\\
Gower Street\\ London\\ WC1E 6BT\\ UK}
\email{Leonid@math.ucl.ac.uk}
\end{addresses}

\begin{abstract}
Counting lattice points inside a ball of large radius in Euclidean space is a classical problem in analytic number theory, dating back to Gauss. We propose a variation on this problem: studying the asymptotics of the measure of an integer lattice of affine planes inside a ball. The first term is the volume of the ball; we study the size of the remainder term. While the classical problem is equivalent to counting eigenvalues of the Laplace operator on the torus, our variation corresponds to the integrated density of states of the Laplace operator on the product of a torus with Euclidean space. The asymptotics we obtain are then used to compute the density of states of the magnetic Schr\"odinger operator.
\end{abstract}

\maketitle

\section{Introduction and Main results}

The first problem we are considering in this paper has several equivalent formulations. 

\subsection{Number theoretic formulation}For $\rho>0$ and $\bk \in \mathbb R^d$, let $B(\rho;\bk)$ be the ball of radius $\rho$ centered at $\bk$. Let $S(\rho;\bk)$ be the number of integer points inside the disk $B(\rho,\bk)\subset \mathbb R^2$. The classical Gauss Circle Problem consists in estimating the remainder term
\bee
\tilde R(\rho;0) = S(\rho;0) - \pi \rho^2
\ene
Hardy and (Edmund) Landau have found lower bounds for this problem, while the current best upper bound is given by Huxley in \cite{Hux}. This problem has also been studied for balls of dimension higher than two, see e.g. \cite{G}, and it is well-known that averaging over the radius of the ball improves regularity of the remainder.

In this paper, we consider a variation on this problem: we estimate the measure of the intersection of affine planes sitting on integer coordinates with balls of large radius in $\mathbb R^d$. More precisely, put 
\bee
A_k:=\mathbb Z^k\times \R^{d-k}\subset\R^d
\ene
and let $B^d(\rho, \bk)$ be a ball in $\mathbb R^d$ of radius $\rho$ centred at $\bk := (\bk_1,\bk_2) \in \mathbb R^k \times \mathbb R^l$, where $k+l = d$. Denote by $S(\rho;\bk_1;d,k)$ the $l$-dimensional volume of the set $B^d(\rho,\bk_1) \cap A_k $. A simple observation shows that we have
\begin{equation} \label{def1}
S(\rho;\bk_1;d,k) = \omega_{l} \sum_{\substack{\gamma \in \Z^k \\ |\gamma - \bk_1| < \rho}}  (\rho^2 - |\gamma - \bk_1|^2)^{l/2}, \end{equation}
where $\omega_{d}$ is the volume of the unit ball in $\mathbb R^d$. One can see that the integral of $\tilde R(\rho,\bk)$ over $\bk_2 \in \mathbb T^l = \mathbb R^l/\mathbb Z^l$, is the same as the remainder term
\begin{equation}
R:=S(\rho;\bk_1;d,k)-\omega_d \rho^d,
\end{equation}
obtained from Equation \eqref{def1}. Our aim is to compute an estimate of $R$ for large values of $\rho$. Before discussing the results, we would like to describe different formulations of this problem. 

\subsection{First spectral theoretic formulation}
Let 
\bee
H=-\Delta+V
\ene
be a Schr\"odinger operator acting in $\R^d$ with a smooth real-valued periodic potential $V$; for simplicity we assume that the lattice of periods $\Gamma=(2\pi\Z)^d$, with dual lattice $\Gamma^\dagger = \mathbb Z^d$ 
Denote the integrated density of states (IDS) of $H$ by $N(\lambda):=N(\la;H)$. It can be defined by the formula
\bee \label{eq:ids}
 N(\la;H) := \lim_{L \to \infty} \frac{\tilde N(\lambda;H_L)}{L^d},
\ene
where $H_L$ is the restriction of $H$ to the cube $[0,L]^d$ with appropriate self-adjoint boundary conditions and $\tilde N(\lambda,H_L)$ is the counting functions of the (discrete) eigenvalues of $H_L$. Note that this parameter $\la$ is related to the parameter $\rho$ of the previous section by $\rho = \sqrt \la$. While this formulation of the IDS is important for Theorem \ref{thm:parland}, for periodic $V$ we use an useful equivalent definition.  

Following \cite{ReedSimon}, we express $H$ as a direct integral 
\bee\label{eq:decper}
H=\int^{\oplus}_{\T^d} H(\bk)\de\bk,
\ene
 Then, one can express $N(\la;H)$ in terms of the counting functions of the fibre operators $H(\bk)$:  
\bee\label{eq:densityperiodic}
N(\lambda):=\frac{1}{(2\pi)^d}\int_{\T^d}N(\la;H(\bk))\de \bk,
\ene
where $N(\la,H(\bk))$ is the eigenvalue counting function of $H(\bk)$. Remarkably, despite the fact that the asymptotic behaviour of $N(\la,H(\bk))$ for fixed $\bk$ and $\la\to\infty$ is very irregular (so that even the precise size of the remainder 
\bee
R(\la;\bk):=N(\la,H(\bk))-C_d\la^{d/2}
\ene
is unknown), integration over all quasimomenta  $\bk\in\T^d:=\R^d/\Z^d$ makes things extremely regular, so that there exists a complete asymptotic expansion of $N(\la)$ in powers  of $\la$ as $\la\to\infty$, \cite{ParnovskiShterenberg1,ParnovskiShterenberg2}.
Here, we have denoted 
\bee
C_d=\frac{\omega_d}{(2\pi)^d} \ \text{and} \ \omega_d=\frac{\pi^{d/2}}{\Gamma(1+d/2)}
\ene
is the volume of the unit ball in $\R^d$. The question we want to study is what would happen if, instead of integrating against all quasimomenta, we integrate over a subset of them, say over an affine plane. We write $\bk=(\bk_1,\bk_2)$, where $\bk_1\in\T^k$, $\bk_2\in\T^l$ and define the partial density of states (PDS) as 
\bee
N_p(\la;\bk_1)=N_p(\la;\bk_1;d,k):=\frac{1}{(2\pi)^d}\int_{\T^l}N(\la,H(\bk_1,\bk_2))\de\bk_2.
\ene

Our aim is to investigate the asymptotic behaviour of the PDS as $\la\to\infty$. Obviously, the regularity at infinity will be improving as $l$ increases and so the larger $l$ is, the more asymptotic terms we are likely to obtain. This asymptotic problem can be treated in two steps: 

Step 1. Obtain the asymptotic behaviour of the PDS for unperturbed operator $H^0:=-\Delta$. More precisely, we want to obtain as good an estimate on 
\bee
R^0(\la;\bk_1;d,k):=N_p^0(\la;\bk_1;d,k)-C_d\la^{d/2}
\ene
as possible (of course, superscript $0$ refers to the fact that we are dealing with the case $V=0$). A simple calculation shows that if $k = 0$, then $R^0(\la;\bk_1;d,0)=0$, so this step is trivial when dealing with the IDS. In the case of $k > 0$ this step becomes quite non-trivial and interesting. Once we have performed this step, we can move to 

Step 2. Compute (or estimate) the difference 
\bee
N_p(\la;\bk_1;d,k)-N_p^0(\la;\bk_1;d,k)
\ene
and try to obtain as many asymptotic terms of it as possible. It follows from a simple computation that
\bee
N_p^0(\la;\bk_1;d,k)=(2\pi)^{-d} S(\sqrt{\la};\bk_1;d,k),
\ene
hence the main aim of this paper deals with the first step of this programme; we intend to perform the second step in a separate publication. 

\subsection{Second spectral theoretic formulation}

Consider the operator $\tilde H=-\Delta+\tilde V$ 
acting on $\T^l\times\R^k$ with a smooth potential $\tilde V: \T^l\times\R^k\to \R$. 
We assume that, as a function on $\R^k$, $\tilde V$ is periodic with the lattice of periods $(2\pi\Z)^k$. Then, we have, from the definition of both the IDS and the PDS that 
\bee
\begin{aligned}
N(\lambda;\tilde H) &= \frac{1}{(2\pi)^k}\int_{\T^l}N(\la;H(\bk_2))\de \bk_2
&= (2\pi)^l N_p(\lambda;0;d,k),
\end{aligned}
\ene
that is to say that the integrated density of states equals the partial density of states up to a constant. If we consider a more general (but also less natural) operator $\tilde H_{\bk_1}$, the domain of which consists of functions on $\T^l\times\R^k$ which become periodic after multiplication by $e^{i\bk_1\bx_1}$, then the IDS of  $\tilde H_{\bk_1}$ equals, again up to the same constant,  $N_p(\la;\bk_1;d,k)$. We would also like to mention that expression \eqref{def1} appears in the study of integer points in anisotropically expanding domains. This has applications in the study of the asymptotic behaviour of the eigenvalue of the Laplace operator on the torus in the adiabatic limit, and was developed in \cite{KY}.

\subsection{Main results}

Our first main result is as follows: 
\begin{thm} \label{thm:asym}
 The error term $R(\rho;\bk_1;d,k)$ satisfies the asymptotic estimates
\begin{equation}
R(\rho;\bk_1;d,k) = \begin{cases}
	   O\left(\rho^{(d-1)/2}\right) & \text{if } k < (d+1)/2, \\
        O\left(\rho^{(d-1)/2}\log \rho\right) & \text{if } k = (d+1)/2, \\
        O\left(\rho^{d-\frac{2k}{1-d+2k}}\right)& \text{if } k > (d+1)/2
       \end{cases}
\end{equation}
uniformly in $\bk_1$.

\end{thm}

\ber
Recall that $R(\rho;\bk_1;d,0)=0$ for all values of $\rho,\bk_1,d$. 
\enr

We do not pretend that all of these estimates are optimal, but some of them are, as can be seen from the following result: 
\begin{thm}\label{lower}
 For $k > 1$ and $\rho$ sufficiently large, there exists a positive constant $C_{d,k}$ and $\bk_1 \in \T^k$ such that
 \begin{equation} 
  R(\rho;\bk_1;d,k) \ge \begin{cases}
			  C_{d,k} \rho^{\frac{d - 1 - \epsilon}{2}} & \text{if } d \equiv 1 \mod 4 \\
                          C_{d,k} \rho^{\frac{d-1}{2}} & \text{else,}
                         \end{cases}
 \end{equation}
where $\epsilon > 0$ is arbitrary. When $d \not \equiv 1 \mod 4$, the lower bound $R(\rho;\bk_1;d,k) \ge C_{d,k} \rho^{\frac{d-1}{2}}$ holds for $k = 1$.
\end{thm}
In particular, this theorem means that for  $1 \le k < \frac{d+1}{2}$ and $d\not \equiv 1 \mod 4$, we cannot get improvements on the upper bounds found in Theorem \ref{thm:asym}. It also means that for $d \equiv 1 \mod 4$, $k \ne 1$, we cannot get  improvements in the exponent. 

\ber
It seems interesting that, after we have integrated $N(\lambda;H(\bk))$ $(d-1)/2$ times, additional integrations do not improve the remainder estimate, until we perform the last ($d$-th) integration, which makes the remainder equal zero. 
\enr

\begin{open}
 The results in \cite{G} imply that for $k = d$, our upper bound is not optimal, but as $d \to \infty$, our upper bound converges to the optimal one, in the sense that $d - (d-\frac{2k}{1-d+2k}) \to 2$. Hence we may ask what is the optimal upper bound for $k \ge \frac{d+1}{2}$.
\end{open}

\subsection{Operators with constant magnetic field}

Another type of problems we consider in this paper is the asymptotic behaviour of the density of states of the (Lev) Landau Hamiltonian (Schr\"odinger operator with constant magnetic field). 

Let $D_j = -i\frac{\D}{\D x_j}$. Then we define the Landau Hamiltonian $H_d$ as the operator acting in $\R^d$ whose action is given by: 
\[
 H_d = (D_1 + x_2)^2 + D_2^2 + \dots + D_d^2.
\]
Of course, only operators $H_2$ and $H_3$ make real physical sense, but for the sake of completeness we will deal with all dimensions.

Let $\Omega^d(\rho)$ for $d\ge 2$ be the parabolic domain in $\mathbb R^d$ given by
\begin{equation}
 \Omega^d(\rho) := \left\lbrace (x_0,x) \in \mathbb R^d : 0 \le x_0 \le \rho -|x|^2 \right\rbrace.
\end{equation}

Defining $P(\rho;d,k)$ analogously to $S(\rho;0;d,k)$, that is,
\begin{equation}
P(\rho;d,k) =  \operatorname{Vol}_{l}(\Omega^d(\rho)\cap A_k),
\end{equation}

one can see that
\begin{equation} \label{eq:blpar}
 P(\rho;d,k) = \sum_{j=0}^{\lfloor \rho \rfloor} S((\rho - j)^{1/2};0;d-1,k-1).
\end{equation}

The IDS $N(\lambda;H_d)$ is related to $P(\rho;d,k)$ by the following proposition. 
\begin{prop} \label{thm:parland}
Let $H_d$ be the $d$-dimensional Landau Hamiltonian. Then, its integrated density of states is given by
\begin{equation}
 N(\la;H_d) = 2^{\frac{-d} 2}\pi^{1-d} P\left(\frac{\la-1}{2};d-1,1\right)
\end{equation}
for $\rho \ge 1$, and $0$ otherwise. 
\end{prop}

We get an asymptotic expression for $P(\rho;d,k)$, via the next theorem. Defining $E_0(\rho) := E_0(\rho,d) = \frac{2}{d+1} \rho^{(d+1)/2} + \frac{1}{2}\rho^{(d-1)/2}$ and
\[
 E_n(\rho):= E_n(\rho,d) = E_0 + \sum_{k=1}^n \frac{B_{2k}}{(2k!)} \frac{\Gamma(\frac{d+1}{2})}{\Gamma(\frac{d+3 - 4k}{2})} \rho^{\frac{d + 1 - 4k}{2}},
\]
we obtain the following theorem.
\begin{thm} \label{thm:para}
 As $\rho\to\infty$, $P(\rho;d,k)$ admits the asymptotic expansions:
 \begin{equation}
  P(\rho;d,1) = \omega_{d-1}E_{\lfloor \frac{d+1}{4} \rfloor}(\rho) + O(1),
 \end{equation}
\begin{equation}
 P(\rho;d,d) = \frac{2\omega_{d-1}}{d+1} + O\left(\rho^{\frac{d^2 - d + 2}{2d}}\right).
\end{equation}
If $k > \frac{d+2}{2}$, we have
\begin{equation}
 P(\rho;d,k) = E_{\lfloor \frac{k-1}{4k-2} \rfloor}(\rho) + O(\rho^{\frac{1}{2}\left(d - 1 - \frac{2k-2}{2k-d}\right)}).
\end{equation}
Finally, if $k \le \frac{d+2}{2}$, 
\begin{equation}
 P(\rho;d,k) =  E_{\lfloor \frac{d-4}{8} \rfloor} (\rho) + O(\rho^{\frac{d+4}{4}}(\log \rho)^{\delta}),
\end{equation}
where $\delta = 1$ if $k = \frac{d+2}{2}$ and $0$ otherwise.
\end{thm}

Replacing the result in Proposition \ref{thm:parland} with the asymptotics in Theorem \ref{thm:para}, we immediately deduce the following corollary.
\begin{cor} \label{thm:corollar}
The integrated density of states of the Landau Hamiltonian on $\mathbb R^3$ admits the asymptotic expansion
\[
 N(\lambda; H_3) = \frac{1}{6\pi^2}\lambda^{3/2} + O(1)
\]
for large enough $\lambda$.
\end{cor}

The rest of the paper is organised as follows: in Section 2 we formulate several  results which will be used in the proof of the main theorems, but we will postpone their proofs until Section 6. In Section 3 we prove the upper bounds in the Laplace case, and in Section 4 we obtain lower bounds. Finally, in Section 5 we deal with the magnetic case. 

\subsection*{Acknowledgments}
The research of J.L. is part of his doctoral studies at Universit\'e de Montr\'eal, under the supervision of Iosif Polterovich.
We are grateful to Zeev Rudnick for outlining the proofs of Lemmas \ref{thm:sand} and \ref{thm:four} in the case $d = 3, k = 2$. We are also grateful to Guillaume Poliquin for providing a generalisation of Lemma \ref{thm:sand} to arbitrary dimension, and for fruitful discussions. We also want to thank Yuri Kordyukov for reading the preliminary version of our manuscript and making useful suggestions as well as for bringing \cite{KY} to our attention.

The research of J.L. was partially supported by the NSERC CGS-M scholarship. The research of L.P. was partially supported by the EPSRC grant EP/J016829/1.

\section{Auxiliary results}

In order to prove Theorem \ref{thm:asym}, it will be useful to give an alternate expression for $S(\rho;\bk_1;d,k)$.  Let us define the function $\chi : \mathbb R^k \to \mathbb R$ as

\begin{equation} \chi(x) = \begin{cases}
              (1-|x|^2)^{l/2}& \text{ if } |x|� < 1, \\
	      0	& \text{otherwise.}
             \end{cases}
\end{equation}
We can then observe that
\begin{equation}
   S(\rho;0;d,k) =  \omega_{l} \rho^{l} \sum_{n \in \mathbb Z^k} \chi(n/\rho).
\end{equation}
 We would like to use Poisson's summation formula
\begin{equation} \label{eq:poisson}
\sum_{n\in \mathbb Z^k} f(n) = \sum_{m \in \mathbb Z^k} \hat f(m)
\end{equation}
with $f = \chi$. This will allow us to get upper bounds for all $\bk_1 \in \T^k$, from the relation
\begin{equation} \label{eq:trans}
 \mathcal F (f(x-\bk_1)) = e^{-2\pi i \bk_1 \cdot \xi} (\mathcal F f)(\xi),
\end{equation}
where $\mathcal F$ is the Fourier transform operator. For the rest of this section, we therefore consider $\bk_1 = 0$, and it will be seen in the proof of Lemma \ref{thm:four} that this assumption is made without loss of generality. In order for Equation \eqref{eq:poisson} to hold we need to smooth out $\chi$. To do so, we will consider its convolution with Friederichs' mollifier $\Psi_\epsilon$. Hence, setting $\chi_\epsilon = \Psi_\epsilon * \chi$ we get that
\begin{equation}
\hat \chi_\epsilon(\xi) = \hat \Psi_\epsilon (\xi) \hat \chi(\xi).
\end{equation}

Theorem \ref{thm:asym} follows from two lemmas. The first one finds asymptotic upper and lower bounds for $S$:
\begin{lem} \label{thm:sand}
 Let $\chi^+_\epsilon$ and $\chi^-_\epsilon$ be defined on $\mathbb R^k$ by

\begin{equation} 
 \chi_\epsilon^\pm(x) = \frac{1}{(1\mp \epsilon)^{l}}\chi_\epsilon((1\mp \epsilon)x).\end{equation}

Then, we have that 
\begin{equation}\label{eq:sand}
\chi^-_\epsilon(x) \le \chi(x) \le \chi^+_\epsilon(x) 
\end{equation}
for all $x\in \mathbb R^k$. Immediately, if we define

\begin{equation}S^\pm_\epsilon(\rho) = \omega_{l}\sum_{n \in \mathbb Z^k} \chi^\pm_\epsilon(n/\rho),\end{equation}

we get that

\begin{equation} S^-_\epsilon(\rho) \le S(\rho) \le S^+_\epsilon(\rho). \end{equation}
\end{lem}
Since $\chi^\pm_\epsilon$ are smooth functions, we can use Poisson's summation formula to compute the asymptotic expansion of $S_\epsilon^\pm$. The second lemma therefore gives the asymptotic expansion of $\hat \chi(\xi)$.
\begin{lem} \label{thm:four}
The Fourier transform of $\chi$ satisfies
\begin{equation} \hat \chi(\xi) = \frac{C}{|\xi|^{(d+1)/2}}\cos \left(2\pi |\xi| - \frac{(d+1)\pi}{4}\right) + O(|\xi|^{-(d+3)/2}) \end{equation}
for some $C>0$ as $|\xi| \to \infty$. Furthermore, its derivative satisfies
\begin{equation}
 \frac{\de}{\de |\xi|} \hat \chi(\xi) = \frac{\tilde C}{|\xi|^{(d+1)/2}}\sin \left(2\pi |\xi| - \frac{d\pi}{4}\right) + O(|\xi|^{-(d+3)/2})
\end{equation}
In particular, the asymptotic behaviour of both $\hat \chi (\xi)$ and its derivative does not depend on the co-dimension $k$. 
\end{lem}

We will postpone the proof of these lemmas until Section \ref{sec:lem}.

\section{Proof of Theorem \ref{thm:asym}} \label{sec:thm}

In this section, we prove Theorem \ref{thm:asym} using both Lemmas \ref{thm:four} and \ref{thm:sand}. We have that
\begin{equation}
 S_\epsilon^-(\rho) \le S(\rho) \le S_\epsilon^+(\rho).
\end{equation}
Let us therefore find asymptotic expansions on $S_\epsilon^\pm$. We shall split those computations in two cases : whether $k \ge (d+1)/2$ or $k < (d+1)/2$

\subsection{Case 1}

Here, we assume that $k\ge (d+1)/2$.
Let us find asymptotic expansions on $S_{\epsilon}^{\pm}$. Since $\chi_\epsilon$ is a smooth compactly supported function of $x$, we may use Poisson's summation formula \eqref{eq:poisson} to obtain
\begin{equation}
 S_\epsilon^\pm = \omega_{d} \rho^{l} \sum_{n \in \mathbb Z^k} \chi^\pm_\epsilon(n/\rho) 
 = \omega_{d} \rho^{d} \sum_{m \in \mathbb Z^k} \hat\chi^\pm_\epsilon(\rho m).
\end{equation}
Since we have that
\begin{equation}
 \hat \chi^\pm_\epsilon(m \rho) = \frac{1}{(1\mp\epsilon)^d}\hat\Psi(\epsilon m \rho)\hat\chi(\frac{m\rho}{1\mp\epsilon}),
\end{equation}
we get, assuming $\epsilon \ll 1/\rho$, that
\begin{equation}
 S_\epsilon^\pm = \omega_{d}\sum_{m\in \mathbb Z^k}(1+O(\epsilon))\rho^d\hat \Psi(\epsilon m \rho)\hat\chi(m\rho) 
  + \bigo{\sum_{m \in \mathbb Z^k} \epsilon m \rho^{d+1}\Psi(\epsilon m \rho) \left|\hat \chi'(m \rho)\right|},
\end{equation}
which directly implies
\begin{equation} \label{eq:Seps}
\begin{aligned}
 S_\epsilon^\pm = \omega_d \rho^d + O(\epsilon \rho^d) &+
 O\left(\sum_{\substack{m\in \mathbb Z^k\\|m|\ne 0}}\rho^d\hat \Psi(\epsilon m \rho)\left|\hat\chi(m\rho)\right|\right)\\ &\qquad + \bigo{\sum_{m \in \mathbb Z^k} \epsilon m \rho^{d+1} \Psi(\epsilon m \rho) \left|\hat \chi'(m \rho)\right|}.
 \end{aligned}
\end{equation}

Observe that $\hat \Psi(\xi) = O(|\xi|^{q})$ for any $q$ whenever $|\xi| >1$ and bounded for $|\xi| \le 1$. Recall from Lemma \ref{thm:four}
that $\hat \chi(\xi) = O(|\xi|^{-(d+1)/2})$. Hence, choosing $q = \frac{d-2k-1}{2}$, the third summand in \eqref{eq:Seps} can be split into two terms, becoming
\begin{equation}
 O\left(\rho^{(d-1)/2}\left[\sum_{\substack{m\in \mathbb Z^k\\1\le |m| \le 1/\epsilon \rho}}\frac{1}{|m|^{(d+1)/2}}+
\sum_{\substack{m\in \mathbb Z^k\\|m| ­­> 1/\epsilon \rho}}\frac{1}{(\epsilon \rho)^{(2k+1-d)/2} |m|^{k+1}}\right]\right).
\end{equation}

The first sum can be estimated by
\begin{equation}
\begin{aligned}
\sum_{\substack{m\in \mathbb Z^k\\1\le |m|\le 1/\epsilon \rho}} \frac{1}{|m|^{(d+1)/2}} 
&\sim \int_1^{1/\epsilon \rho}\frac{r^{k-1}}{r^{(d+1)/2}}\de r\\ 
&= \begin{cases} 
O\left(\left(\epsilon \rho\right)^{\frac{d+1-2k}{2}}\right) & \text{if } k > (d+1)/2, \\
O(\log \epsilon \rho) & \text{if } k = (d+1)/2 .
\end{cases}
\end{aligned}
\end{equation}

The second sum can be estimated by
\begin{equation}
\begin{split}
 & \sum_{\substack{m\in \mathbb Z^k\\ |m|\ge 1/\epsilon \rho}} \frac{1}{(\epsilon \rho)^{(2k-d+1)/2}|m|^{k+1}}\\
& \sim \int_{1/\epsilon \rho}^{\infty}\frac{1}{(\epsilon \rho)^{(2k-d+1)/2}}\frac{r^{k-1}}{r^{k+1}}\de r = 
O\left((\epsilon \rho)^{\frac{d+1-2k}{2}}\right).
\end{split}
\end{equation}

As for the last summand, it is easy to see with the same computations and using $\hat \Psi(\xi) = \bigo{|\xi|^{\frac{d - 2k - 3}{2}}}$ that the extra power of $\epsilon\rho$ exactly compensates the extra power of $m$, and we have that the asymptotic behavior in $\epsilon \rho$ is the same for all for summands whenever $k> (d+1)/2$.
Furthermore, when equality holds, the polynomial component is the same. Therefore, we have to choose $\epsilon = \rho^{-j}$ such that
 \begin{equation}
\epsilon \rho^d = \rho^{(d-1)/2}(\epsilon \rho)^{\frac{d+1-2k}{2}}.
 \end{equation}

This is achieved exactly when
\begin{equation}
 j = \frac{2k}{1 - d  + 2k}.
\end{equation}

This gives us the announced asymptotic estimates when $k\ge (d+1)/2$, that is
\begin{equation}
S(\rho) = \begin{cases}\omega_d \rho^d + O(\rho^{d-\frac{2k}{1-d+2k}}) & \text{if } k > (d+1)/2, \\
        \omega_d \rho^d + O(\rho^{\frac{d-1}{2}}\log \rho) & \text{if } k = (d+1).
       \end{cases}
\end{equation}

\subsection{Case 2}

We now assume that $k< (d+1)/2$. 
In this case, we have that the sum converges with $\hat \Psi = O(1)$. Hence, the asymptotic expansion for $S_\epsilon^\pm$
simplifies to

\begin{equation}
 S_\epsilon^\pm = \omega_d \rho^d + O(\epsilon \rho^d) +
 O\left(\rho^{(d-1)/2}\sum_{\substack{m\in \mathbb Z^k\\|m|\ne 0}}\frac{1}{|m|^{(d+1)/2}}\right) + \bigo{\rho^{(d-1)/2}\sum_{\substack{m \in \mathbb Z^k \\ m \ne 0}}\frac{\epsilon \rho \hat\Psi(\epsilon m \rho)}{|m|^{(d-1)/2}}}.
\end{equation}
The third sum converges and the last one as well if $k < \frac{d-1}{2}$. In that case, choosing $\epsilon = \rho^{-(d+1)/2}$ satisfies Theorem \ref{thm:asym}, and choosing $\epsilon$ smaller does not improve the estimate. If $k = \frac d 2$ or $k = \frac{d-1}{2}$, using $\hat \Psi(\xi) = \bigo{|\xi|^{-1}}$ for $m > (\epsilon\rho)^{-1}$ yields the same result, finishing the proof.

Note that Equation \eqref{eq:trans} ensures that these estimates hold for all $\bk_1 \in \T^k$.

\section{Lower bounds}

%\subsection{Lower bounds in $\bk_1$}

Let us first follow the argument given in \cite{DahlbergTrubowitz1982} for $d = k = 2$. The beginning of the argument is the same, which we add for completeness. Since $R(\rho;\bk_1)$ is periodic in $\bk_1$ with respect to $\Gamma$, we can compute its Fourier coefficients, obtaining
\begin{equation}
 \begin{aligned}
&\int_{\T^k} R(\rho ; \bk_1) e^{-2\pi i\bk_1 \cdot \gamma} \de \bk_1 = \int_{\T^k}\left(- \omega_d \rho^d + \rho^{l}\sum_{\gamma \in \Gamma} \chi\left(\frac{\gamma - \bk_1}{\rho}\right)e^{-2\pi i\bk_1 \cdot \gamma}\right) \de \bk_1 \\
&= \int_{\R^k} \rho^l \chi\left(\frac{\bk_1}{\rho}\right) e^{-2\pi i\bk_1 \cdot \gamma}\de \bk_1 \\
&= \rho^d\left[\frac{C}{(\rho |\gamma|)^{(d+1)/2}}\cos \left(2\pi\rho |\gamma| - \frac{(d+1)\pi}{4}\right) +  O(|\rho\gamma|^{-(d+3)/2})\right],
 \end{aligned}
\end{equation}
from Lemma \ref{thm:four}. Additionally, we have that
\begin{equation}
 \int_{\T^k} R(\rho;\bk_1) \de \bk_1 = 0.
\end{equation}
Hence, for all $\gamma \in \Gamma \setminus \lbrace 0 \rbrace$, we have that
\begin{equation} \label{eq:max}
 \begin{aligned}
&  \int_{\T^k} |R(\rho;\bk_1)| \de \bk_1\\ &\ge \max \left(\left|\int_{\T^k} R(\rho;\bk_1)e^{-2\pi i \bk_1 \cdot \gamma} \de \bk_1 \right|,\left|\int_{\T^k} R(\rho;\bk_1)e^{-4\pi i \bk_1 \cdot \gamma} \de \bk_1 \right| \right) \\
  &\ge C \frac{\rho^{\frac{d-1}{2}}}{\gamma^{\frac{d+1}{2}}} \max \left(\left|\cos\left(2\pi \rho |\gamma| - \frac{(d+1)\pi}{4}\right)\right|,\frac{1}{2^{­\frac{d+1}{2}}}\left|\cos\left(4\pi \rho |\gamma| - \frac{(d+1)\pi}{4}\right)\right| \right) \\
  & \qquad - c \frac{\rho^{\frac{d-1}{2}}}{\gamma^{\frac{d-1}{2}}} \\
 \end{aligned}
\end{equation}
for $C,c$ positive constants whose value can change throughout. Whenever $d \not \equiv 1 \mod 4$, we have that
\begin{equation}
 0 < \inf_{x \in \R} \max\left(\left| \cos\left(x - \frac{(d+1) \pi}{4}\right)\right|,\left| \cos\left(2x - \frac{(d+1) \pi}{4}\right)\right|\right),
\end{equation}
hence in that case, fixing $\gamma \in \Gamma$, we conclude that there exists $r^*$ such that for all $r\ge r^*$
\begin{equation}
  \int_{\T^k} |R(\rho;\bk_1)| \de \bk_1 \ge C \rho^{\frac{d-1}{2}}.
\end{equation}
We conclude that whenever $d \not \equiv 1 \mod 4$,
\begin{equation}
 \sup_{\bk_1 \in \T^k} R(\rho;\bk_1) \ge C \rho^{\frac{d-1}{2}}.
\end{equation}

The remaining case, that is when $d \equiv 1 \mod 4$ is more subtle. We will use results found in \cite{ParnovskiSobolev2001}[Theorem 3.1, Lemma 3.3]. Indeed, from Equation \eqref{eq:max}, we have
\begin{equation}
  \int_{\T^k} |R(\rho;\bk_1)| \de \bk_1 \ge C \frac{\rho^{\frac{d-1}{2}}}{\gamma^{\frac{d+1}{2}}} \left|\cos\left(2\pi \rho |\gamma| - \frac{\pi}{2}\right)\right| - c \frac{\rho^{\frac{d-1}{2}}}{\gamma^{\frac{d-1}{2}}}.
\end{equation}
From Lemma 3.3 in \cite{ParnovskiSobolev2001}, we know that, if $k\ge 2$, for all $\epsilon > 0$, there exists $\rho_0>0$ and $\alpha \in (0,1/2)$ such that for all $\rho > \rho_0$ there exists $\gamma \in \Gamma$ such that $|\gamma| < (2\pi\rho)^\epsilon$ and the distance from $2\rho\gamma$ to an integer is greater than $\alpha$. Choosing such a $\gamma$ bounds $\cos(2\pi\rho|\gamma| - \pi/2)$ away from $0$, and we get that
\begin{equation}
 \int_{\T^k} |R(\rho;\bk_1)| \de \bk_1 \ge C \rho^{\frac{d-2}{2} - \epsilon\frac{d+1}{2}}.
\end{equation}
Since $\epsilon > 0$ is arbitrary, we get the desired result.

\section{An application to the Landau Hamiltonian}
\subsection{The Landau Hamiltonian} \label{sec:equiv}

Decomposing $H_d = H_2 \oplus D_{d-2}$, we can first study the problem
\[
 H_2 u = \lambda u.
\]
Consider the definition \eqref{eq:ids} for  $N(\lambda;H_d)$, with periodic boundary conditions for $x_1$ and Dirichlet boundary conditions for $x = (x_2,\dotsc,x_d)$.

For $H_2$, we can write the solutions as $u(x_1,x_2) = e^{\frac{2\pi i n}{L}x_1}f(x_2)$, which reduces the problem to solving the eigenvalue problem
\[
 \left((\xi_1 + x_2)^2 + D_2^2\right) f(x_2) = \lambda f(x_2).
\]
This is a shifted quantum harmonic oscillator. We have that $\sigma(H_2) = \lbrace 2j+1 : j \in \mathbb N\rbrace$, each with infinite multiplicity. It is a standard computation, see e.g. \cite{Nakamura}, that
\bee
N(\lambda;H_2) = \frac{1}{2\pi}\left\lfloor \frac{\lambda - 1 }{2} \right \rfloor,
\ene
for $\lambda \ge 1$, and $0$ otherwise.
Extending the methods of \cite{Nakamura} to higher dimensions, it is again a simple computation to show that for $\lambda \ge 1$,
\bee
 N(\lambda;H_d) = \frac{\omega_{d-2}}{(2\pi)^{d-1}}\sum_{n=0}^{\lfloor \frac{ \lambda - 1 }{2} \rfloor} (\lambda-2n-1)^{(d-2)/2}.
\ene

Thus, from the definition of $P(\rho;d,k)$, we have indeed that 
\bee
 N(\lambda;H_d) = 2^{\frac{-d} 2}\pi^{1-d} P\left(\frac{\la-1}{2};d-1,1\right).
\ene

\subsection{Computations for general paraboloids}

In this section we prove Theorem \ref{thm:para}. Consider the expression
\begin{equation} \label{eq:paraa}
 P(\rho;d,k) = \sum_{j=1}^{\lfloor \rho \rfloor} S((\rho-j)^{1/2};0;d-1,k-1).
\end{equation}

By Theorem \ref{thm:asym}, we have  
\[
 \sum_{j=0}^{\lfloor \rho \rfloor} S((\rho - j)^{1/2};0;d-1,k-1) = \sum_{j=0}^{\lfloor \rho \rfloor} \left(\omega_{d-1} (\rho - j)^{(d-1)/2} + O(X(\rho))\right),
\]
where
\begin{equation} \label{eq:X}
 X(\rho) = \begin{cases}
      \rho^{\frac{1}{2}\left(d-1 - \frac{2k-2}{2k-d}\right)} & \text{if } k> (d+2)/2, \\
      \rho^{(d-2)/4} \log \rho & \text{if } k = (d+2)/2, \\
      \rho^{(d-2)/4}  & \text{if } 1 < k < (d+2)/2, \\
      0 & \text{if } k=1.
     \end{cases}
\end{equation}
Comparing with the integral, we get that for all $X$ as defined above,
\begin{equation}
  \sum_{j=0}^{\lfloor \rho \rfloor} X(\rho) = O(\rho X(\rho)).
\end{equation}

For any $d$, we can use the Euler-Maclaurin formula : 
\begin{equation}
\begin{aligned}
& \sum_{n=a}^b f(n) = \int_a^b f(x) \de x + \frac{f(a) + f(b)}{2} \\
 &\quad + \sum_{k=1}^p \frac{B_{2k}}{(2k)!} 
\left(\frac{\de^{2k-1} f}{\de x }\bigg|_{x = b} -\frac{\de^{2k-1} f}{\de x }\bigg|_{x = a}\right) + O\left(\int_a^b \left|\frac{\de^{2p} f}{\de x^{2p} }\right|_{x=t} \de t\right),
\end{aligned}
\end{equation}
for any integer $p\ge 1$, where $B_k$ is the $k$th Bernoulli number. Note that for integer $a$, 
\begin{equation}
 \sum_{j=0}^a (a-j)^{(d-1)/2} = \sum_{j=0}^a j^{(d-1)/2}.
\end{equation}
Hence, by the Euler-Maclaurin formula, we get that
\begin{equation}
\begin{aligned}
& \sum_{j=0}^a (a-j)^{(d-1)/2} \\&= \int_0^a t^{(d-1)/2} \de t + \frac{a^{(d-1)/2}}{2} + \sum_{k\le \frac{d+1}{4}} \frac{B_{2k}}{(2k!)} \frac{\Gamma(\frac{d+1}{2})}{\Gamma(\frac{d+3 - 4k}{2})} a^{\frac{d + 1 - 4k}{2}}  + O(a^{-1/2}).
 \end{aligned}
\end{equation}
 Obviously, when $d$ is odd, this last sum is actually finite and the error term $0$.

When $\rho$ is not an integer, we can write $\rho = a + \tau$, where $\tau$ is the fractional part. In that case, using the Euler-Maclaurin formula again, we get
\begin{align*}
& \sum_{j=0}^a  (a + \tau - j)^{(d-1)/2} = \sum_{j=0}^a (j + \tau)^{(d-1)/2}\\
&= \int_0^a (t + \tau)^{(d-1)/2} \de t + \frac{1}{2} \left(\tau^{(d-1)/2} + \rho^{(d-1)/2}\right) \\
& \quad+ \sum_{k\le \frac{d+1}{4}} \frac{B_{2k}}{(2k!)} \frac{\Gamma(\frac{d+1}{2})}{\Gamma(\frac{d+3 - 4k}{2})} \left(\rho^{\frac{d + 1 - 4k}{2}} - \tau^{\frac{d + 1 - 4k}{2}}\right) + O(\tau)\\
 &= \frac{2}{d+1}\left(\rho^{(d+1)/2} - \tau^{(d+1)/2}\right) +  \frac{1}{2} \left(\tau^{(d-1)/2} + \rho^{(d-1)/2}\right) \\
&\quad+ \sum_{k\le \frac{d+1}{4}} \frac{B_{2k}}{(2k!)} \frac{\Gamma(\frac{d+1}{2})}{\Gamma(\frac{d+3 - 4k}{2})} \left(\rho^{\frac{d + 1 - 4k}{2}} - \tau^{\frac{d + 1 - 4k}{2}}  \right)+ O(\tau) .
\end{align*}

Let us observe that 
\begin{align*}
& \lim_{\rho \to \infty} \frac{\frac{-2}{d+1}\rho^{(d+1)/2} + \sum_{j=0}^a  (\rho - j)^{(d-1)/2}}{\frac{1}{2} \rho^{(d-1)/2}} \\&=\lim_{\rho \to \infty} \frac{-\frac{4}{d+1}\tau^{(d+1)/2} + \tau^{(d-1)/2} + \rho^{(d-1)/2} + O(\rho^{(d-3)/2})}{\rho^{d-1}} \\
&= 1.
\end{align*}
This is because $\tau = O(1)$. Similarly, if we define $E_0 = \frac{2}{d+1} \rho^{(d+1)/2} + \frac{1}{2}\rho^{(d-1)/2}$ and
\[
 E_n = E_0 + \sum_{k=1}^n \frac{B_{2k}}{(2k!)} \frac{\Gamma(\frac{d+1}{2})}{\Gamma(\frac{d+3 - 4k}{2})} \rho^{\frac{d + 1 - 4k}{2}},
\]
we get that
\begin{align*}
 &\lim_{\rho \to \infty} \frac{-E_n + \sum_{j=0}^a  (\rho - j)^{(d-1)/2}}{\rho^{(d-1)/2 - 2n - 1}} = \frac{B_{2(n+1)}}{(2(n+1))!} (\frac{d-1}{2})_{2n+1}
\end{align*}
 whenever $(d-1)/2 - 2n -1 ­> 0$, after which point the contribution of the fractional remainder $\tau$ gets more important than the denominator. Hence, we obtain the asymptotic expansion
\begin{equation}\label{eq:rem}
\begin{aligned}
 \sum_{j=0}^{\lfloor \rho \rfloor}  (\rho - j)^{(d-1)/2} &= \frac{2}{d+1} \rho^{(d+1)/2} + \frac{1}{2}\rho^{(d-1)/2} \\
&\quad + \sum_{k \le \frac{d-3}{4}} \frac{B_{2k}}{(2k!)} \frac{\Gamma(\frac{d+1}{2})}{\Gamma(\frac{d+3 - 4k}{2})} \rho^{\frac{d + 1 - 4k}{2}} + O(\tau).
\end{aligned}
\end{equation}

When $k = 1$, we already have that $X(\rho) = 0$. Therefore, we have that 
\begin{align*}
 \frac{P(\rho,d,1)}{\omega_{d-1}} &= \frac{2}{d+1} \rho^{(d+1)/2} + \frac{1}{2}\rho^{(d-1)/2} \\
&\quad + \sum_{1 \le k < \frac{d-3}{4}} \frac{B_{2k}}{(2k!)} \frac{\Gamma(\frac{d+1}{2})}{\Gamma(\frac{d+3 - 4k}{2})} \rho^{\frac{d + 1 - 4k}{2}} + O(\tau),
\end{align*}
from which we recover a (quite sharp) asymptotic integrated density of states for the magnetic Hamiltonian $H_{d+1}$.

Let us combine  equations \eqref{eq:X} and \eqref{eq:rem}. When $k = d$, we get that the error term from $X$ is greater than $\frac{d-1}{2}$, and as such,
\[
 P(\rho;d,d) = \frac{2\omega_{d-1}}{d+1} + O\left(\rho^{\frac{d^2 - d + 2}{2d}}\right).
\]
When $k > \frac{d+2}{2}$, we get that
\begin{align*}
 \frac{P(\rho;d,k) }{\omega_{d-1}}&= \frac{2}{d+1} \rho^{\frac{d+1}{2}} + \frac{1}{2}\rho^{\frac{d-1}{2}} \\ 
 & \quad + \sum_{1 \le j < \frac{k-1}{4k-2}} \frac{B_{2j}}{(2j)!} \frac{\Gamma(\frac{d+1}{2})}{\Gamma(\frac{d+3 - 4j}{2})} \rho^{\frac{d + 1 - 4j}{2}} + O(\rho^{\frac{1}{2}\left(d - 1 - \frac{2k-2}{2k-d}\right)}).
\end{align*}

Finally, when $k \le \frac{d+2}{2}$, we get that 
\begin{align*}
 \frac{P(\rho;d,k)}{\omega_{d-1}}  &= \frac{2\omega_{d-1}}{d+1} \rho^{\frac{d+1}{2}} + \frac{1}{2}\rho^{\frac{d-1}{2}} \\ 
 & \quad + \sum_{1 \le j < \frac{d+4}{8}} \frac{B_{2j}}{(2j)!} \frac{\Gamma(\frac{d+1}{2})}{\Gamma(\frac{d+3 - 4j}{2})} \rho^{\frac{d + 1 - 4j}{2}} + O(\rho^{\frac{d+4}{4}}(\log \rho)^{\delta}),
\end{align*}
where $\delta = 1$ if $k = \frac{d+2}{2}$ and $0$ otherwise.

\section{Proofs of auxiliary results} \label{sec:lem}

\subsection{Smoothing of the cut-off function}

Let us define a smooth, even bump function $\psi$ in $C_c^\infty(\mathbb R)$, supported in $[-1,1]$, such that
the integral
\begin{equation} \int_0^\infty \psi(r) r^{k-1} \de r = \frac{1}{V_{k-1}}, \end{equation}
where $V_{k-1}$ is the area of the unit sphere in $\mathbb R^k$.

Using this function, we can define the radial bump function $\Psi_\epsilon$ on $\mathbb R^k$, of total mass $1$ to be given by
 \begin{equation} \Psi_\epsilon (x) = \frac{1}{\epsilon^{k}} \psi (\frac{|x|}{\epsilon}). \end{equation}

Let $\Psi := \Psi_1$ and $\chi_\epsilon(x) = \Psi_\epsilon (x) * \chi(x)$. Its Fourier transform is given by

\begin{equation} \hat{\chi_\epsilon}(\xi) = \hat\Psi(\epsilon \xi) \hat \chi (\xi). \end{equation}

 Let $\chi^+_\epsilon$ and $\chi^-_\epsilon$ be defined on $\mathbb R^k$ by

\begin{equation} 
 \chi_\epsilon^\pm(x) = \frac{1}{(1\mp \epsilon)^{l}}\chi_\epsilon((1\mp \epsilon)x).\end{equation}

We can now proceed with the proof of Lemma \ref{thm:sand}.

\bep
To show that $\chi_\epsilon^- (x) \le \chi(x) \le \chi_\epsilon^+$, the idea is to obtain $\chi_\epsilon^\pm (x)$ by averaging $\chi(x)$ on a ball of radius $0< \epsilon < x$ about each $x$. To do so, first notice that
\begin{equation}
\begin{aligned}
 \chi_\epsilon(x) &\le \sup_{|t| \le \epsilon } (\chi(x-t))  \int_{\mathbb R^k} \Psi_\epsilon (x) \de x \\
 &=\begin{cases}
     1  &\text{if } |x| \le \epsilon,\\
     (1-(|x|-\epsilon)^2)^{\frac{l}{2}} & \text{if } \epsilon \le |x| \le 1 +  \epsilon.
   \end{cases}
\end{aligned}
\end{equation}
If we show that
\begin{equation}\label{eq:key} 
\chi_\epsilon(x) \le (1+\epsilon)^{l}\chi\left(\frac{x}{1+\epsilon}\right),
\end{equation}
 we get the desired lower bound. Indeed, taking $y = \frac{x}{1+\epsilon}$ in the preceding
equation yields 
\begin{equation}
 \chi(y) \ge \frac{1}{(1+\epsilon)^{l}}\chi_\epsilon((1+\epsilon)y) = \chi_\epsilon^-(y).
\end{equation}
Therefore, it only remains to show that \eqref{eq:key} holds for all $x\in \mathbb R^k$. First note that if $|x| \ge 1+\epsilon$, both sides are $0$.
We shall split the remaining cases in $|x| \le \epsilon$ and $\epsilon < |x| < 1+\epsilon$. 

 Restricting ourselves to the first case, if $|x| = \epsilon$, we get that
\begin{equation}
 \begin{aligned}
  (1+\epsilon)^{l} \chi\left(\frac{x}{1+\epsilon}\right) &= (1+\epsilon)^{l}\left(1-\frac{\epsilon^2}{(1+\epsilon)^2}\right)^{\frac{l}{2}} \\ 
&= (1+2\epsilon)^{\frac{l}{2}} \\
&\ge 1 \\
&\ge \chi_\epsilon^{-}(x)
 \end{aligned}
\end{equation}
Since $\chi(\frac{x}{1+\epsilon})$ is a decreasing function of $|x|$, we conclude that \eqref{eq:key} holds for $0\le |x| \le \epsilon$. 

In the case where $\epsilon < |x| \le 1+\epsilon$, we need to show that
\begin{equation}
(1-(|x|-\epsilon)^2)^{\frac{l}{2}} \le (1+ \epsilon)^{l}\left(1-\frac{|x|^2}{(1+\epsilon)^2}\right)^{\frac{l}{2}}.
\end{equation}
It is equivalent to show that $1-(|x|-\epsilon)^2 \le (1+\epsilon)^2 - |x|^2$. This is the case if
\begin{equation}
 \begin{aligned}
  1 - |x|^2 + 2|x|\epsilon - \epsilon^2 &\le 1 + 2\epsilon + \epsilon^2 - |x|^2 \\ 
\Leftrightarrow 2|x|\epsilon &\le 2\epsilon(1+\epsilon)\\
\Leftrightarrow |x| &\le 1+\epsilon.
 \end{aligned}
\end{equation}
Since the last line is true by hypothesis, we can conclude that the left-hand side inequality of \eqref{eq:sand} is true. 

In order to get an upper bound on $\chi(x)$, we proceed in a similar fashion, averaging $\chi_\epsilon(x)$ on a ball of radius $\epsilon$ around $x$, which yields
\begin{equation}
 \begin{aligned}
  \chi_\epsilon(x) &\ge \inf_{|t|<\epsilon}\chi(x-t) \\
&\ge \begin{cases}
      \left(1-(|x|+\epsilon)^2\right)^{\frac{l}{2}} & \text{if } |x| < 1-\epsilon ,\\
      0 & \text{otherwise.}
     \end{cases}
 \end{aligned}
\end{equation}
As we did before, it suffices to show that
\begin{equation}
 \chi_\epsilon(x) \ge (1-\epsilon)^{l}\chi\left(\frac{x}{1-\epsilon}\right).
\end{equation}
Notice that the left hand side of that equation is $0$ whenever $|x|\ge1-\epsilon$. Like before, we see that
\begin{equation}
(1-(|x|+\epsilon)^2)^{\frac{l}{2}} \ge (1-\epsilon)^{l}\left[1-\left(\frac{|x|}{1-\epsilon}\right)^2\right]^{\frac{l}{2}}
\end{equation}
is equivalent to $|x|<1-\epsilon$. This concludes the proof.
\enp

\subsection{Fourier transform of $\chi$}

\bep
Let us compute $\hat \chi(\xi)$. We will split the cases $k=1$, $k=2$, and $k>2$. If $k=1$, then

\begin{equation}
\begin{aligned}
 \hat \chi (\xi) &= \int_{-1}^1 (1-x^2)^{(d-1)/2}e^{-i2\pi x\xi} \de x \\
 &= \frac{C}{|\xi|^{d/2}}J_{d/2}(2\pi|\xi|) \\
 &= \frac{C}{|\xi|^{(d+1)/2}}\cos \left(2\pi|\xi| - \frac{(d+1)\pi}{4} \right) + O(|\xi|^{(d+3)/2}),
\end{aligned}
\end{equation}
 using \cite{GradshteynRyzhik1994}[Eq.3.387  and 8.451], which is the desired result.

We also obtain that, following \cite{GradshteynRyzhik1994}[Eq. 3.621]

\begin{equation}
 \hat \chi (0) = 2^d B(\frac{d+1}{2},\frac{d+1}{2}).
\end{equation}
Using identities of the Gamma function, we get that
\begin{equation}
\omega_{l }2^d B(\frac{d+1}{2},\frac{d+1}{2}) = \frac{\pi^{d/2}}{\Gamma(\frac{d}{2}+1)} = \omega_d, 
\end{equation}
which is the desired value.

If $k=2$, then the Fourier transform is given by

\begin{equation} \hat \chi (\xi) = \int_{\mathbb R^2} \chi(x) e^{-i2\pi x\cdot \xi} \de x. \end{equation}

Working in polar coordinates, we get that

\begin{equation}
\begin{aligned}
\hat \chi (\xi) &= \int_0^1 \int_0^{2\pi} r(1-r^2)^{(d-2)/2} e^{-i2\pi r |\xi| \cos \theta} \de \theta \de r \\
&=\int_0^1 r(1-r^2)^{(d-2)/2} J_0(2\pi|\xi|r) \de r \\
&= \frac{C}{|\xi|^{d/2}}J_{d/2}(2\pi |\xi|) \\
&= \frac{C}{|\xi|^{(d+1)/2}}\cos \left(2\pi|\xi| - \frac{(d+1)\pi}{4}\right) + O(|\xi|^{(d+3)/2}),
\end{aligned}
\end{equation}
which is the desired result. \cite{GradshteynRyzhik1994}[Eq. 8.411, 6.567 and 8.451] were used respectively for an integral formula for the Bessel function, its integral, and its asymptotic expansion.

We also obtain that 
\begin{equation}
 \hat \chi (0) = \frac{2\pi}{d}.
\end{equation}
Using identities of the Gamma function, we get that
\begin{equation}
 \omega_{l}\frac{2\pi}{d} = \frac{\pi^{d/2}}{\Gamma(\frac{d}{2}+1)} = \omega_d, 
\end{equation}
which is the desired value.
Finally, if $k>2$, then, working in spherical coordinates, we get that the Fourier transform of $\chi$ is, for some constant $C$,

\begin{equation}
 \begin{aligned}
  \hat \chi (\xi) &= C \int_0^1 \int_0^{\pi} r^{k-1}(1-r^2)^{l/2} \sin^{k-2}\theta e^{-i2\pi r |\xi| \cos \theta} \de \theta \de r \\
  &= \frac{C}{|\xi|^{(k-2)/2}}\int_0^1 r^{k/2}(1-r^2)^{l/2} J_{(k-2)/2}(2\pi|\xi|r) \de r \\
  &= \frac{C}{|\xi|^{(k-2)/2}}\frac{1}{|\xi|^{(l+2)/2}}J_{d/2}(2\pi|\xi|) \\
  &= \frac{C}{|\xi|^{(d+1)/2}}\cos \left(2\pi|\xi| - \frac{(d+1)\pi}{4}\right) + O(|\xi|^{(d+3)/2}).
 \end{aligned}
\end{equation}

using \cite{GradshteynRyzhik1994}[Eq. 8.411] in the first line, which is the desired result.

Additionnally, we have that

\begin{equation}
\begin{aligned}
 \hat \chi (0) &= \operatorname{Vol}(S^{k-1})\int_0^1 r^{k-1} (1-r^2)^{(d-k)/2} \de r \\
&= \frac{\pi^{k/2}B(\frac{k}{2},\frac{d-k+2}{2})}{\Gamma(\frac{k}{2})}.
\end{aligned}
\end{equation}
Using identities of the Gamma function, we get that
\begin{equation}
 \hat \chi(0) \omega_{d-k} = \omega_d
\end{equation}
which is once again the desired value.

One can note that in each of those cases, we ignored the trigonometric term to get an upper bound, considering it to be $1$. Hence, since translation by $\bk_1$ is simply multiplication by a complex exponential in Equation \eqref{eq:poisson}, it can be ignored in just the same fashion.

Finally, we get the result for the derivative using the identity $J_{\nu}' = \frac 1 2 (J_{\nu - 1} - J_{\nu + 1})$ and basic trigonometric identities. This completes the proof of Lemma \ref{thm:four}.
\enp

\end{document}